# Acceptable area of optimal control for a multidimensional system

S N Masaev[1,4], G A Dorrer[1,2] and V V Cyganov[3]

[1] Siberian Federal University, pr. Svobodnyj, 79, Krasnoyarsk, 660041, Russia  
[2] Department of System Analysis and Operations Research, Office L-409, 410, 31, Krasnoyarsky Rabochy Av., Reshetnev Siberian State University of Science and Technology, Krasnoyarsk, 660037, Russia  
[3] Laboratory 57, V. A. Trapeznikov Institute of Control Sciences of RAS, 65 Profsoyuznaya street, Moscow, 117997, Russia

[4] E-mail: faberi@list.ru

**Abstract**. A research shows the phase space of the system. The phase space of such a system is determined by the development structure of four subsystems with different objective functions. The control loop of such a system is formed. Using the control loop, optimal control is generated. The dynamic control region is calculated on the basis of a matrix determining the structure of the development of a multidimensional dynamic system. It was established that the optimal distribution by R. Bellman's method allows increasing the increase in the value of the objective function over 5 years by 85% from the initial value with a decrease in the distributed resource by 20%. It is believed that the construction of railways leads to an increase in gross regional product by 2%, but the authors proved that it is possible to increase this figure to 7-8%

## 1. Introduction

The formation of an optimal solution is one of the main issues in control theory. Despite the fact that the optimal control problem and its solutions are widely known, difficulties arise in calculating the optimal solution in multidimensional dynamic non-stationary systems [1-8]. Optimal control in such systems is not obvious, because of their size and complexity of the structure of the interaction of parameters in time. It is especially important to evaluate the optimal control for various target functions in the system for the implementation of short-term projects. The paper investigates the problem of the optimal allocation of resources in the development of special economic zones and its impact on the return on railway construction for the subject of the Russian Federation [9,10]. In an applied sense, the calculation of optimal control is a rather complicated task. The control object has a large dimension and with each time step changes its structure under the influence of the implementation of its own goals and the influence of environmental factors. The purpose: to calculate the allowable area of control of a multidimensional system.

## 2. Method

The multidimensional dynamical system is given by the equation

$$x(t+1) = A(t)x(t) + B(x(t)u(t)) + v(t) \qquad (1)$$

$$y(t)=G((x(t),u(t),v(t))$$
$$x(t)\in R^n, u(t)\in R^m, v(t)\in R^l, y(t)\in R^k$$

where $x(t)\in R^n$ – vector of state signals of the object, $u(t)\in R^m$ – vector control actions at the system level, $v(t)\in R^l$ – environmental disturbance vector, $y(t)\in R^k$ – observation vector $y(t)=Hx(t)$, $u=W(q)y$. We find the transfer matrix $W$ by solving a problem of the form: $J=J(W)\to \inf_{w\in\Omega^*}$. $Q$ – forward shift operator, $H=[h_{ij}]$ – $K\times N$ – matrix defining system monitoring structure, $J$ – quality parameter for system states at various $q$, determining the admissible set of transfer matrices $\Omega^*$ taking into account all the requirements for the dynamics of the process. Requirements are determined through matrices: $A=[a_{ij}]$ – $N\times N$ is the matrix that determines the speed of development and the structure of the system, $B=[b_{ij}]$ – $N\times M$ is the matrix that defines the structure of the control action on the system.

## 3. Practical task

The object of research is the optimal control of the special economic zone (SEZ) of the Krasnoyarsk Territory, located according to plans on the adjacent territory of Yemelyanovo Airport, 40 km. from the city of Krasnoyarsk. The SEZ is the general territory where enterprises are located and carry out their economic activities (residents of the SEZ) [11]. The special economic zone is the territory where enterprises (residences) operate with preferential tax regimes. Many domestic and foreign works have been devoted to the characterization of activities and control problems of the SEZ [12-34]. The head of the constituent entity of the Russian Federation manages the SEZ through the preferential regime of tax rates provided by subsidies and determines who to include or exclude from the number of residents. According to the development strategy of the Krasnoyarsk Territory until 2020, it is planned to create a cargo hub on the basis of the Emelyanovo airport in Krasnoyarsk [9,10] and a special economic zone. According to studies of a large foreign consulting company, the promising development sectors of the Krasnoyarsk Territory and the SEZ include: the woodworking industry (1st resident) [35], the baby food (dry mix) industry (2nd resident), petrochemicals (3rd - resident) and radio electronics (4th resident). Subsystems are residents of SEZ. The data structure, economic activities and ongoing projects of residents for modeling are described for a period of 5 years: the structure of the project and owners, the characteristics of the products being produced, market analysis, organizational plan, marketing plan, production plan, resource analysis, environmental impact assessment, analysis project risks, financial model [35-37]. This project involves the construction and use of the railway line for the pass-fat train to the airport and back - 40 km. The railway line is also used for the carriage of goods by residents of the SEZ. The project parameters contain the estimated cost of building various route options, the cost of buying land in accordance with the cadastral plan, the number of passengers transported, the nomenclature and weight of the goods transported, the need for rail transport in the surrounding settlements. Control of the SEZ is carried out by the head of the subject of the Russian Federation.

The activity of the special economic zone is determined by equation 1. The interpretation of the parameters of equation (1) is preserved. $X=[X^1,X^2,X^3,X^4]$ - a special economic zone (SEZ) of the constituent entity of the Russian Federation ($X$) consists of residents of the SEZ: a woodworking enterprise, an oil refinery, a radio engineering enterprise, and manufacturer of dry infant formula. The SEZ is controlled through $u(t)=(u_1(t),u_2(t),...u_m(t))^T$ is the $M$ vector, where $u_i(t)$ is the above-mentioned control actions by the state in the form of changes in interest rates of taxes, subsidies to the industry and benefits aimed at increasing the gross regional product ($y$-GRP) at the time $t$, $u_1$ - income tax rate to the budget of the subject of the Russian Federation, $u_2$ - income tax rate in the consolidated budget of the Russian Federation, $u_{12}=u_1+u_2$, $u_{12}$ - general income tax rate, $u_3$ - transport tax, $u_4$ - property tax, $u_5$ - electricity cost, $u_6$ - measures to reduce logistics costs for the movement of inventory items, reimbursement of the cost of transportation of products (subsidies), $u_7$ - the average value of contributions to the Pension Fund of the Russian Federation (PF RF), the Territorial Fund for Compulsory Medical Insurance (MHIF), the Federal Fund for Compulsory Medical Insurance (FFOMS), the Social Insurance Fund of the Russian Federation (FSS RF), $u_8$ - the price of renting a land plot (land rent), $u_9$ – the price of renting a plot with resources. The work of the SEZ over 5 years is estimated by the gross regional product $y(t)$ 130 billion rubles, $u_6$=16 billion rubles. subsidies for transportation costs (construction of a railway line from the SEZ on the territory of the Krasnoyarsk airport to the city of Krasnoyarsk 40 km long) [9,10], $u_{12}$=1.92 billion rubles. A figure 1 shows the effect of all tax rates $u_n$ on y SEZ activities is disclosed in detail.

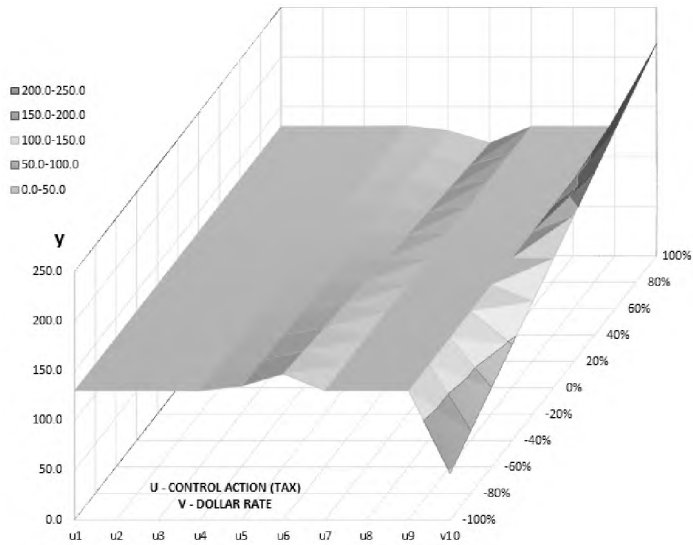

**Figure 1.** Change ($y_i$) SEZ from the investment policy regimes of the subject of the Russian Federation (Krasnoyarsk Territory).

It is required to select the control $u_{12}=(u_{12}^1;u_{12}^2;u_{12}^3;u_{12}^4)$ in such a way as to maximize the function $J$. $q$=5 years, taking into account the construction of the railway from the city of Krasnoyarsk to the SEZ. Transfer matrices $W$ are determined by the individual parameters of the activity of each resident of the SEZ, depending on the size of the income tax rate (benefits) $u_{12}$. The calculation was carried out for 4.6 million values in the author's software package [38-40]. It is required to select the control $u_{12}=(u_{12}^1;u_{12}^2;u_{12}^3;u_{12}^4)$ in such a way as to maximize the function J. Figure 2 and table 1 show the optimal distribution of subsidies is found by R. Bellman's method [41].

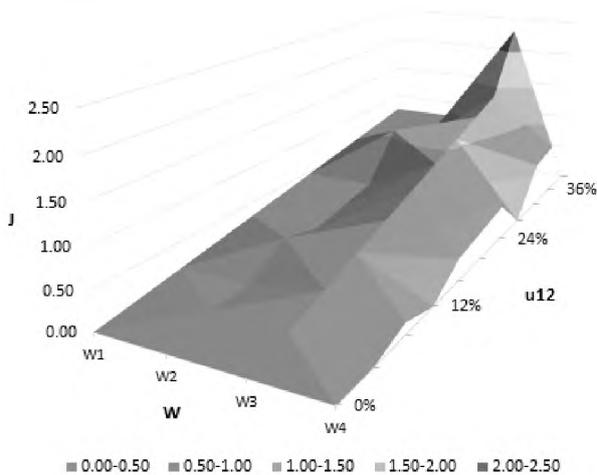

**Figure 2.** Optimal control area (multiple transfer matrices $W$).

Table 1. The values of the objective functions depending on the optimal distribution of subsidies ($q$=5 years).

| $u_{12}^k$ | $W(q)$ | $y$ | $J$ |
|---|---|---|---|
| 0% | 0 | 149.8 | 0 |
| 4% | 0.38 | 142.7 | 0.5 |
| 8% | 0.77 | 150.0 | 1.0 |
| 12% | 1.15 | 142.7 | 1.5 |
| 16% | 1.54 | 150.0 | 1.8 |
| 20% | 1.92 | 150.0 | 2.2 |
| 24% | 2.30 | 150.0 | 2.6 |
| 28% | 2.69 | 142.7 | 3.4 |
| 32% | 3.07 | 150.0 | 3.4 |
| 36% | 3.46 | 150.0 | 3.9 |
| 40% | 3.84 | 143.7 | 5.6 |

According to the law of Russia, the income tax rate (benefits) $u_{12}$ cannot exceed 20%. It can be seen from the calculation (see table 1) that it is enough to reduce $u_{12}$ by 4% to 16% so that the railway takes its maximum estimation, taking into account the construction of the railway. With the condition, at $u_{12}$=16%, the conditional maximum of the target functions $J$=1.8, and $y$=150 billion rubles are achieved with the condition (optimal control of $J$) that the first, second, fourth residents subsidized 0.31 billion rubles, the third resident subsidized 0.61 billion rubles. Total is 1.54 billion rubles (0.31 billion rubles, 0.31 billion rubles, 0.61 billion rubles, 0.31 billion rubles). Given the above limitations, stimulating the activities of residents of the SEZ, it is possible to increase profits by 0.24 billion rubles from 1.54 billion rubles up to 1.8 billion rubles. Then the assessment of the effectiveness of the control function of the SEZ will be 0.25 billion rubles.

## 4. Conclusion
A model of a dynamic system is formalized. The control loop of such a system is defined. To control the dynamic system, the R. Belman method is applied. The optimal control area of the dynamic system is calculated. The calculation was carried out for 4.6 million values in the author's software package [38-40]. Modeling showed optimal control of the SEZ $u_{12}$ - 1.54 billion rubles=(0.31; 0.31; 0.61; 0.31). Over 5 years, an increase in indicators of the Krasnoyarsk Territory from the construction of the railway line (infrastructure project) will amount to: $y$ - GRP from 130 to 150 billion rubles, $J$ - profit from 1.54 to 1.8 billion rubles. The payback of the infrastructure project in the SEZ will be reduced from 8 years to 3 months. The research goal achieved. The work was performed as a continuation of a comprehensive research of the optimal control of multidimensional economic facilities. [42-46].